\documentclass{elsart}

\usepackage{amssymb}
\journal{ArXiv}

\newtheorem{theorem}{Theorem}[section]
\newtheorem{lemma}[theorem]{Lemma}

\theoremstyle{definition}

\newtheorem{example}[theorem]{Example}

\newcommand*{\pd}[2]{\frac{\partial #1}{\partial #2}}

\begin{document}
\begin{frontmatter}




\title{From reflections to elliptic growth}
\author[Tanya]{Tatiana  Savina }
\ead[url]{https://people.ohio.edu/savin/}
\address[Tanya]{Department of Mathematics,  Ohio University,  USA}

\begin{abstract}
{We use reflections involving analytic Dirichlet and Neumann data on a real-analytic curve in order  to find a representation of solutions to  Cauchy problems for harmonic functions in the plane. 
We apply this representation for finding solutions to  Hele-Shaw problems. 
We also generalize the results 
by deriving the corresponding formulae for Helmholtz equation and applying them to elliptic growth.
}
\end{abstract}

\begin{date}
\date{}
\end{date}

\end{frontmatter}
\maketitle

\section{  Introduction  }

A connection between Hele-Shaw problem and analytic continuation has been mentioned in literature. The aim of this work is to show this connection explicitly. Our starting point was the celebrated Schwarz reflection principle.

Reflections are usually associated with continuation of boundary value problems. They were extensively studied by many researchers (see for example  \cite{ek}, \cite{khavinsons},
\cite{dima}, \cite{nonloc}, \cite{shapiro} and references therein). 
Analytic continuation of solutions to Cauchy problem for elliptic equations was studied in \cite{Mi1}.
In this paper, we use reflections for harmonic  functions as well as for solutions to the Helmholtz equation satisfying to nonhomogeneous Dirichlet or Neumann boundary conditions to obtain continuation for a Cauchy problem and applying it to Hele-Shaw flows.

The mathematical formulation of the two-phase Hele-Shaw problem is as follows.
Let $\Omega _2 (t) \subset {\mathbb R}^2$ with a boundary $\Gamma (t)$ at time $t$ be a 
simply-connected bounded domain occupied by a fluid with a constant viscosity $\nu _2$,   and let $\Omega _1 (t)$ be the region  ${\mathbb R}^2\setminus {\bar \Omega}_2(t)$ occupied by a different fluid of viscosity $\nu _1$.
Consider the two-phase Hele-Shaw flow forced by sinks and sources:
\begin{equation}\label{1}
{\bf v}_j=-k_j\nabla p_j, \qquad j=1,2, 
\end{equation}
where the pressure $p_j$ is a harmonic function almost everywhere in the region $\Omega _j (t)$, satisfying boundary conditions
\begin{eqnarray}\label{2}
p_1(x,y,t)=p_2(x,y,t) \quad \mbox{on} \quad \Gamma (t),\\
-k_1\pd{p_1}{n}=-k_2\pd{p_2}{n} =v_n \quad \mbox{on} \quad \Gamma(t).\label{3}
\end{eqnarray}
Here ${\bf v}_j$ is a velocity vector of fluid $j$, $k_j=h^2/12\nu _j$, and $h$ is the gap width of the Hele-Shaw cell.
Equation (\ref{2}) states the continuity of the pressure under the assumption of negligible
 surface tension.  
 Equation (\ref{3}) means that the normal velocity of the boundary itself coincides with the normal velocity of the fluid
 at the boundary.

The free boundary $\Gamma (t)$ moves due to the sources and sinks  
 located in both regions.
One may add a natural physical assumption that the fluid flux generated by the system of sources and sinks is finite. That allows no more than the logarithmic growth of the fluid pressure near a point source/sink or at  infinity.

In a special case when $k_1=k_2=k$, a two-color problem,  some  solutions may be constructed as an analytic continuation of one-phase problem. For instance, consider the following  problem
\begin{eqnarray}\label{2_1}
\Delta p(x,y,t)=\mu \quad \mbox{near} \quad \Gamma (t),\\
p(x,y,t)=\phi \quad \mbox{on} \quad \Gamma (t), \label{2_2}\\
-k\pd{p}{n} =v_n \quad \mbox{on} \quad \Gamma(t).\label{3}
\end{eqnarray}
Note that (i) $\phi =0$ in the case of one-phase problem with a negligible surface tension, (ii) $\phi$ is a known function in the case of one-phase problem, which describes either a surface tension or an external force \cite{EE}, \cite{KMV}, \cite{mcdonald}, and (iii) $\phi$ is unknown in the case of the two-color problem.
However in the special case when the problem (\ref{2_1})-(\ref{3}) is considered as an inverse problem with a prescribed dynamics, this dynamics could occur in a way that the interface belongs to a given family of curves defined by a zero level set of the pressure, then $\phi$ vanishes.

Here we assume that the support of  of the distribution $\mu$, which describes sinks and sources, is  not located near the free boundary. Then one can  consider (\ref{2_1})-(\ref{3}) as a Cauchy problem for a harmonic function.

The structure of the paper is as follows. In Section 2, we review a  reflection formula for a nonhomogeneous  Dirichlet condition, derive the corresponding formula for the Neumann conditions, and discuss the Cauchy problem for the Laplace's equation. 
In section 3, we  show how the derived formula helps to treat Hele-Shaw problems when the flows are induced by sinks and sources as well as by a change in the gap width of a Hele-Shaw cell.
In section 4, we discuss nonhomogeneous reflection formulae for the Dirichlet and Neumann conditions for the Helmholtz equations. Then, we discuss  a representation of solutions to the corresponding Cauchy problem. We show its relationship with elliptic growth 
in section 5, and we formulate the conclusions in section 6.

\section{Reflections and Cauchy problem for Laplace's equation}
\subsection{Reflection for nonhomogeneous Dirichlet data}

Reflection operator is generally an integro-differential operator,
which reduces in the simplest case to the celebrated local (point to
point) Schwarz symmetry principle for harmonic functions, which
can be stated as follows.

Let $\Gamma \subset \mathbb{R}^2$ be a non-singular real analytic
curve and a point $K\in \Gamma$. Then, there exists a neighborhood $U$ of
$K$ and an anti-conformal mapping $R:U\to U$ which is identity on
$\Gamma$, permutes the components $U_1, U_2$ of $U\setminus \Gamma$
and relative to which any harmonic function $u(x, y)$ defined near
$\Gamma $ and vanishing on $\Gamma $ (the homogeneous Dirichlet condition) is
odd,
\begin{equation} \label{E:0.1}
         u(x_0 ,y_0 ) =-u(R(x_0 ,y_0 )).
\end{equation}

In the case of nonhomogeneous Dirichlet data, $u=\phi (x,y)$ on $\Gamma$, when function $\phi$ is holomorphically continuable into $\mathbb{C}^2$ near $\Gamma$, formula (\ref{E:0.1}) involves also values of function $u$ at two more  points located on the  complexification  $\Gamma _{\mathbb{C}}$ of the curve $\Gamma$. All four points then create a so-called Study's rectangle \cite{226}.

To describe the Study's rectangle, consider a complex domain $W$ in the space $\mathbb{C}^2$ to which
the function $f$ defining the curve $\Gamma :=\{f(x,y)=0\}$  can be
 analytically continued such that $W\cap \mathbb{R}^2=U$. Using the
change of variables $z=x+iy$, $w=x-iy,$ the equation of the
complexified curve $\Gamma _{\mathbb{C}}$ can be rewritten in the
form
$$
f\left ( \frac{z+w}{2}, \, \frac{z-w}{2i} \right ) =0,
$$
and if $grad \, f(x,y)\ne 0$ on $\Gamma$, can be also rewritten in
terms of the Schwarz function and its inverse, $w=S(z)$ and
$z=\widetilde{S}(w)$ \cite{davis}. The mapping
$R$ mentioned above can be expressed in terms of the Schwarz function as follows,
\begin{equation} \label{E:0.2}
R(x_0 ,y_0 )=R(z_0 )=\overline{S(z_0 )}.
\end{equation}
Using the above notations, the reflection formula for harmonic functions subject to conditions 
$u_|{_{\Gamma}}=\phi$ can be written as the Study's rectangle:
\begin{equation} \label{E:0.1n}
 u(z_0 ,w_0 ) + u(\widetilde{S}(w_0) ,S(z_0 ))=\phi (\widetilde{S}(w_0),w_0)+\phi (z_0, S(z_0)),
\end{equation}
where $z_0=x_0+iy_0$ and $w_0=x_0-iy_0$.

Indeed, the general solution to the Laplace's equation,
$
\frac{\partial ^2 u}{\partial z\partial w}=0,
$ 
has the form $u(z,w)=g(z)+f(w)$. Thus, 
\begin{eqnarray}\nonumber
 u(z_0 ,w_0 ) +u(\widetilde{S}(w_0) ,S(z_0 ))=
\Bigl [(g(z_0)+f(w_0)\Bigr ]+\Bigl [g (\widetilde{S}(w_0))+f (S(z_0))\Bigr ]=\\
\Bigl [g (\widetilde{S}(w_0))+f(w_0)\Bigr ]+\Bigl [(g(z_0)+f (S(z_0))\Bigr ]=
\phi (\widetilde{S}(w_0),w_0)+\phi (z_0, S(z_0)).\nonumber
\end{eqnarray}

\subsection{Reflection for nonhomogeneous Neumann data}

Let $\Gamma \subset \mathbb{R}^2$ be a non-singular real analytic
curve and $K\in \Gamma$. Then, there exists a neighborhood $U$ of
$K$ and an anti-conformal mapping $R:U\to U$ which is identity on
$\Gamma$, permutes the components $U_1, U_2$ of $U\setminus \Gamma$
and relative to which any harmonic function $u(x, y)$ defined near
$\Gamma $ and satisfying to the homogeneous Neumann condition on $\Gamma$, $\pd{u}{n}=0$, is
even,
\begin{equation} \label{E:0.1b}
         u(x_0 ,y_0 ) =u(R(x_0 ,y_0 )).
\end{equation}

The nonhomogeneous version of formula (\ref{E:0.1b}) is given by the following lemma.
\begin{lemma}\label{L1}
Let $u(x,y)$ be a harmonic function defined near curve $\Gamma$ and satisfying Neumann boundary condition  $\pd{u}{n}=\psi$ on $\Gamma$, where $n$ is a normal and $\psi (x,y)$ is a function that can be holomorphically continued into $\mathbb{C}^2$ near $\Gamma$, then the following reflection formula holds
 \begin{equation} \label{E:0.2n}
         u(x_0 ,y_0 ) = u(R(x_0 ,y_0 ))+i \int\limits _{ \widetilde{S}(w_0)}^{z_0}\psi(z, S(z)) \sqrt{ S^{\prime}}\, dz,
\end{equation}
where $S^{\prime}$ means the derivative of $S$ with respect to $z$.
\end{lemma}
{\it Proof of Lemma \ref{L1}. }
A general solution to the Laplace's equation,
$$
\frac{\partial ^2 u}{\partial z\partial w}=0,
$$ 
has the form $u(z,w)=g(z)+f(w)$. Then
$$
\pd{u}{n}=-\frac{i}{\sqrt{ S^{\prime} }}\Bigl (\pd{g}{z}-\pd{f}{w}S^{\prime}\Bigr )=\psi ,
$$
which implies
$$
\pd{g}{z}=i\psi\sqrt{ S^{\prime} }+\pd{f}{w}S^{\prime}.
$$
Integrating the latter formula along $\Gamma _{\mathbb{C}}$, 
$$
\int\limits _{ (\widetilde{S}(w_0),w_0)}^{(z_0, S(z_0))}\pd{g}{z}\, dz=
i \int\limits _{ (\widetilde{S}(w_0),w_0)}^{(z_0, S(z_0))}\psi(z, S(z)) \sqrt{S^{\prime}}\, dz
+\int\limits _{ (\widetilde{S}(w_0),w_0)}^{(z_0, S(z_0))}\pd{f}{w}S^{\prime}\, dz,
$$
we have
$$
g(z_0)-g(\widetilde{S}(w_0))=
i \int\limits _{ (\widetilde{S}(w_0),w_0)}^{(z_0, S(z_0))}\psi(z, S(z)) \sqrt{S^{\prime}}\, dz
+f(S(z_0))-f(w_0),
$$
which results in
$$
u(z_0,w_0)-u(\widetilde{S}(w_0),S(z_0))=
i \int\limits _{ \widetilde{S}(w_0)}^{z_0}\psi(z, S(z)) \sqrt{S^{\prime}}\, dz
$$
and finishes the proof.

Remark that the sign in front of the integral depends on the choice of the branch of 
$\sqrt{S^{\prime}}$. 
\begin{example}{\rm
Consider a harmonic function satisfying the condition $\pd{u}{n}=\pd{u}{y}=-2y + \alpha$ on the $x$-axis.
Here $\alpha$ is an arbitrary constant. Taking into account that $S(z)=z$ and choosing 
$\sqrt{S^{\prime}}=-1$, formula (\ref{E:0.2n}) reduces to
$$
u(x_0,y_0)=u(x_0,-y_0)+2\alpha\, y_0.
$$
}
\end{example}
\begin{example}{\rm
Consider a harmonic function satisfying the condition $\pd{u}{n}=\pd{u}{r}=\beta$ on the circle
$x^2+y^2=a^2$.
Here $a$ and $\beta$ are arbitrary constants, and $r=\sqrt{x^2+y^2}$. The Schwarz function is
$S(z)=a^2/z$. Choosing $\sqrt{S^{\prime}}=-ia/z$,
Then formula (\ref{E:0.2n}) reduces to
$$
u(x_0,y_0)=u\Bigl( \frac{a^2x_0}{x_0^2+y_0^2},\frac{a^2y_0}{x_0^2+y_0^2}  \Bigr )+
ab\ln{(x_0^2+y_0^2)}-2ab\ln{a}.
$$
}
\end{example}

\subsection{Representation of solutions to Cauchy problem for harmonic functions}

In this subsection, we use the nonhomogeneous reflections mentioned above to obtain a representation of solutions to Cauchy problems.

\begin{theorem}\label{T1}
Solution to the Cauchy problem
\begin{eqnarray}\label{2_1c}
\frac{\partial ^2 u}{\partial z\partial w}=0\quad \mbox{near} \quad \Gamma_{\mathbb{C}},\\
u(z,w)=\phi (z,S(z))\quad \mbox{on} \quad \Gamma_{\mathbb{C}} ,\\
\pd{u}{z}-\pd{u}{w}\pd{S}{z}=i \psi (z, S(z)){\sqrt{ S^{\prime} }}\quad \mbox{on} \quad \Gamma_{\mathbb{C}}\label{3c}
\end{eqnarray}
has a representation
\begin{equation}\label{CR}
u(z_0,w_0)=\frac{ \phi (\widetilde{S}(w_0),w_0)+\phi (z_0, S(z_0))}{2}+
\frac{i}{2} \int\limits _{ \widetilde{S}(w_0)}^{z_0}\psi(z, S(z)) \sqrt{S^{\prime}}\, dz.
\end{equation}
\end{theorem}
{\it Proof of Theorem \ref{T1}. }
The proof immediately follows from (\ref{E:0.1n}) and (\ref{E:0.2n}).

\section{Reformulation of Hele-Shaw problems in terms of free boundary data}

\subsection{Dynamics induced by sinks and sources}

Consider a problem of finding a harmonic function $p(x,y,t)$ defined near free boundary $\Gamma (t) \subset\mathbb{R}^2$ subject to conditions (\ref{2_2})-(\ref{3}). In this case 
$\psi=-v_n/k$. Taking into account $v_n=-i\dot S(z,t)/\sqrt{4S^{\prime}(z,t)}$ \cite{Mi2},
$$
\psi =\frac{i\,\dot S(z,t)}{2k\sqrt{S^{\prime}(z,t)}},
$$
where $\dot S$ means the derivative of $S$ with respect to $t$.
Then formula (\ref{CR}) yields
\begin{equation}\label{p_0}
p(x,y,t)=\frac{ \phi (\widetilde{S}(\bar z,t),\bar z,t)+\phi (z, S(z,t),t)}{2}-
\frac{1}{4k} \int\limits _{ S^{-1}({\bar z},t) }^{z}\dot S (s,t)\, ds.
\end{equation}
The author believes that a version of this formula first appeared in \cite{Mi2}.

In the case when $\phi =0$ (negligible surface tension) one may use formula (\ref{p_0}) for finding a solution to a one-phase Hele-Shaw inverse problem. 

\begin{example}{\rm
 If the interface belongs to a family of circles, $x^2+y^2=a^2(t)$, the computation of the integral in (\ref{p_0}) gives an expression for the pressure
$$
p=-\frac{a\dot a}{2k}\ln\Bigl ((x^2+y^2)/a^2\Bigr ),
$$
which has singularities at zero and infinity. These singularities correspond to sinks and sources.
Remark that in the case of circular dynamics, surface tension does not affect the dynamics. In the presence of surface tension $\phi=\gamma/a(t)$, where $\gamma$ is a constant surface tension coefficient. Thus,  the pressure is
$$
p=-\frac{a\dot a}{2k}\ln\Bigl ((x^2+y^2)/a^2\Bigr )+\frac{\gamma}{a}.
$$
In this example, the interface is moving due to a presence of a point sink/source at zero or/and a source/sink at infinity depending on whether an interior, exterior or a two-color problem is considered.
}
\end{example}
\begin{example}{\rm 
If the interface belongs to a family of ellipses 
\newline
$\Gamma (t)=\left\{ \frac{x^2}{a(t)^2}+\frac{y^2}{b(t)^2}=1\right\}$ with semi-axes $a(t)$ and $b(t)$, where   $a(0)>b(0)$,
the corresponding Schwarz function is
\begin{equation}\label{SE}
S\left( z,t\right)
=\Bigl ( \bigl ( a(t)^2+b(t)^2 \bigr ) z-2a(t)b(t)\sqrt{z^2-d(t)^2}\Bigr )/d(t)^2,
\end{equation}
where $d(t)=\sqrt{a(t)^2-b(t)^2}\,$ is the length of a half of the inter-focal distance.
Consider two different scenario with a negligible surface tension, $\phi =0$.

(i) If the eccentricity of the ellipse does not change with time, the latter implies that the ratio $a(t)/b(t)=const$, the pressure does not have
no more than logarithmic singularities (including the singularity at infinity): 
\begin{equation}
p=
 -\frac{1}{2k_j}\pd{(ab)}{t}\Bigl (\ln |z+\sqrt{z^2-d^2}| -\ln (a+b)\Bigr ),\label{wellm1}
\end{equation}
or
$$
p=
 -\frac{1}{2k_j}\pd{(ab)}{t}\Bigl (\ln \sqrt{(x+\alpha )^2 (1+y^2/\alpha ^2) } -\ln (a+b)\Bigr ),
$$
where 
$$
\alpha ^2=\Bigl ( x^2-y^2-d^2 +\sqrt{(x^2-y^2-d^2)^2+4x^2y^2}  \Bigr ) /2.
$$
The inter-focal distance, $d(t)=2b(t)\sqrt{a^2(0)/b^2(0)-1}$, of such an ellipse changes.
The interface of a two-color problem is moving due to 
 a point sink/source  at infinity and a  source/sink
 distribution with density $\mu _2(x,t)=\frac{2ab}{d^2 k_2} \, \partial _t(\sqrt{d^2-x^2})$
along the inter-focal segment.
 Here $\int _{-d}^{d} k_2 \mu _2(x,t)\, dx= \pi\partial _t(ab)=\dot A$, where $A$ is the area of the ellipse.

(ii) If the area of the elliptical inclusion  does not change in time, that is, $a(t)b(t)=const$,  the pressure $p$ is defined by
$$
p=\frac{1}{k}\Re \Bigl \{-\frac{z^2}{4}\pd{}{t}\Bigl (  \frac{a^2+b^2}{d^2}   \Bigr ) 
+\frac{ab\,z}{2}\sqrt{z^2-d^2}\,\pd{}{t}\Bigl (  \frac{1}{d^2}   \Bigr )
\Bigr \} - \frac{b^2a \dot a}{k_1 d^2},
$$
or
$$
p=\frac{1}{2k}\Bigl \{\frac{(y^2-x^2)}{2}\pd{}{t}\Bigl (  \frac{a^2+b^2}{d^2}   \Bigr ) 
+\frac{ab\,x(\alpha ^2-y^2)}{\alpha}\,\pd{}{t}\Bigl (  \frac{1}{d^2}   \Bigr )
\Bigr \} - \frac{b^2a \dot a}{k_2 d^2}.
$$
The interior flow is generated by the density 
$$
\mu = \frac{ab\,\partial _t (d^2)}{k_2d^4}\,\,\frac{(2x^2-d^2)}{\sqrt{d^2-x^2}},
$$
supported on the inter-focal segment. Such a density changes sign along the inter-focal segment, so the area of the ellipse does not change in time: if $a(t)$ increases with time, the ellipse becomes ``thiner''.
}
\end{example}

\subsection{Dynamics induced by change of the gap width}

Formula (\ref{CR})
can be also applied to the case when the free boundary $\Gamma (t)$ is moving due to change in the gap between the plates in a Hele-Shaw cell instead of the presence of sinks and sources.
A statement of the problem with a time-dependent gap $h(t)$ between the plates
was mentioned in \cite{EEK} among other generalized Hele-Shaw flows. The one-phase (interior)
version of this problem was considered in \cite{tian}, while the two-phase problem in \cite{Avi}.
In the case of a one-phase problem or a two-color problem with a changing gap, the Laplace's equation (\ref{2_1}) is replaced with a Poisson's equation. Thus, the problem for the pressure is
\begin{eqnarray}\label{1g}
\Delta p=\frac{1}{k}\frac{\dot h(t)}{h(t)} \quad \mbox{near} \quad \Gamma (t),\\
p(x,y,t)=0 \quad \mbox{on} \quad \Gamma (t), \label{2g}\\
-k\pd{p}{n} =v_n \quad \mbox{on} \quad \Gamma(t).\label{3g}
\end{eqnarray}
Using the substitution 
\begin{equation}\label{sub}
p(x,y,t)=p_h(x,y,t)+\frac{1}{4k}\frac{\dot h(t)}{h(t)}(x^2+y^2),
\end{equation} 
we obtain the following problem for $p_h(x,y,t)$:
\begin{eqnarray}\label{2_1g}
\Delta p_h(x,y,t)=0 \quad \mbox{near} \quad \Gamma (t),\\
p_h(x,y,t)=-\frac{1}{4k}\frac{\dot h(t)}{h(t)}(x^2+y^2) \quad \mbox{on} \quad \Gamma (t), \label{2_2g}\\
\pd{p_h}{n} =-\frac{1}{k}(v_n+\frac{1}{4}\frac{\dot h(t)}{h(t)}\pd{}{n}(x^2+y^2) \quad \mbox{on} \quad \Gamma(t).\label{3g}
\end{eqnarray}
Thus, 
to find a harmonic function $p_h$, one has to solve the Cauchy problem (\ref{2_1c})-(\ref{3c}) with
$$
\phi (x,y,t)=-\frac{1}{4k}\frac{\dot h(t)}{h(t)}(x^2+y^2), \quad
\psi(x,y,t)=-\frac{1}{k}(v_n+\frac{1}{4}\frac{\dot h(t)}{h(t)}\pd{}{n}(x^2+y^2))
$$
or
\begin{equation}\label{pp}
\phi(z,\bar z,t)=-\frac{1}{4k}\frac{\dot h(t)}{h(t)}z\bar z, \quad
\psi(z,S(z,t),t)=\frac{i}{2k\sqrt{S^{\prime}}}\Bigl (\dot S+\frac{\dot h(t)}{2h(t)}(S -
zS^{\prime})
\Bigl ).
\end{equation}
Remark that taking into account (\ref{pp}), formula (\ref{CR}) can be rewritten as
\begin{eqnarray}\nonumber
& p_h=-\frac{1}{8k}\frac{\dot h(t)}{h(t)}(\bar z \,\widetilde{S}(\bar z,t)+z\, S(z,t))\\
&-\frac{i}{4k}\int\limits _{ \widetilde{S}(\bar z,t)}^{z}\Bigl ( \dot S(s,t)+\frac{\dot h(t)}{2h(t)}(S(s,t) -s\partial _s S(s,t))\Bigr )ds,\nonumber
\end{eqnarray}
which reduces  to
\begin{equation}\label{p_gap}
p_h=-\frac{\dot h(t)}{4kh(t)}\,\bar z \,\widetilde{S}(\bar z,t)
-\frac{1}{4k}\int\limits _{ \widetilde{S}(\bar z,t)}^{z}\Bigl ( \dot S(s,t)+\frac{\dot h(t)}{h(t)}(S(s,t) \Bigr )ds.
\end{equation}
Formula (\ref{p_gap}) along with (\ref{sub}) solves a one-phase and a two-color inverse problems when the interface is moving due to a time-dependent change of the gap width.

\begin{example}{\rm
 When the  interface is circular, $x^2+y^2=a^2(t)$,
it is easy to see that due to volume conservation, $2\dot a+\frac{\dot h a}{h}=0$, the integrand in (\ref{p_gap}) vanishes. Thus,
$$
p_h(x,y,t)
=-\frac{a^2}{4k}\frac{\dot h(t)}{h(t)}.
$$
}
\end{example}
\begin{example}{\rm 
 In the case when the interface belongs to a family of confocal ellipses, 
$d(t)=d(0)=const$, plugging the expression for the Schwarz function (\ref{SE}) into formula
(\ref{p_gap}), we obtain
\begin{equation}
 p_h =\frac{1}{4k_j}\Bigl ( (x^2-y^2)\frac{\dot a\, d_0^2}{a(a^2-d_0^2)}+2\dot a a  
 \Bigr ).
\end{equation}
}
\end{example}

\section{Reflections and Cauchy problem for Helmholtz equation}

\subsection{Reflection for nonhomogeneous Dirichlet data}

First, we review  the reflection formula for solutions $u(x,y)$ to elliptic differential equations with real-analytic coefficients and the Laplacian in the principal,
\begin{equation}\label{geneq}
\hat L u\equiv \Delta  _{x,y}u+a\pd{u}{x}+b\pd{u}{y}+cu=0\mbox{ near
}\Gamma ,
\end{equation}
where
$$
u(x,y)_{\mid _\Gamma }=0 ; \,a,\,b,\, c\,\,\mbox{ are real-analytic functions of }x,\,y.
$$

 Let $\Gamma$ be an algebraic curve, $S(z)$
is  analytic in
$\mathbb{C}$ except for finitely many algebraic singularities.

\begin{theorem} \label{P:0.1} \cite{nonloc}
Under the above assumptions, the following reflection formula holds:

\begin{eqnarray}\label{RF}
\, & u(P) \,   =\,  -\,  c_0(P,\Gamma )\, u(Q)\,  +\\
& \frac{1}{2i}\int\limits   _{\Gamma }^Q \Bigl (
\bigl
\{
u\pd{V^D}{x}-V^D\pd{u}{x}-auV^D \bigr \} dy
      - \bigl \{u\pd{V^D}{y}-V^D\pd{u}{y}-buV^D \bigr \} dx
\Bigr ) .\nonumber
\end{eqnarray}

where $P=(x_P ,y_P)$ and $Q=R(P)$.
Here the integral is independent on the path
joining an arbitrary point on $\Gamma$ with the point $R(x_0 ,y_0)$.
\begin{eqnarray}\nonumber
c_0(P,\Gamma ) &=\frac{1}{2}\Bigl \{ \exp \Bigl [\int\limits _{z_P} ^{{\widetilde
S}(w _P)} B(t,S(z_P))dt+\int\limits _{w _P} ^{
S(z _P)} A(z_P ,\tau)d\tau \Bigr ] \\
& +  \exp \Bigl [\int\limits _{w _P} ^{
S(z _P)} A({\widetilde S}(w _P),\tau )d\tau+\int\limits _{z _P} ^{
{\widetilde S}(w  _P)} B(t ,w _P)dt \Bigr ] \Bigr \} ,
\end{eqnarray}
where
$$
A(z,w )=\frac{1}{4}\bigl [a(x,y)+ib(x,y) \bigr ],
\quad
B(z,w )=\frac{1}{4}\bigl [a(x,y)-ib(x,y) \bigr ],
$$
$$
C(z,w )=\frac{1}{4}c(x,y), \quad
V^D=V^D(x_P ,y_P ,x,y)=V_1^D (x_P ,y_P ,x,y)-V_2^D (x_P ,y_P ,x,y).
$$
\end{theorem}

Functions $V_j$ are solutions of the Cauchy-Goursat problems:
\begin{eqnarray}\nonumber
 {\hat L}^* V_j^D=0, \qquad j=1,2\, ,\qquad\\
 {V_j^D}_ {|_{\Gamma _{\mathbb{C}}}}={\mathfrak R}_{|_{\Gamma _{\mathbb{C}}}},  \qquad  j=1,2\, ,\qquad \\
 V_1^D= \exp \Bigl \{\int\limits _{w _P} ^{
w} A(\widetilde{S}(w ),\tau )d\tau + 
\int\limits _{z _P} ^{z} B(t,w  )dt\Bigr \} \,\,
 \mbox{on } \widetilde{l}_1=\{\widetilde{S}(w )=z_P\},\\
 V_2^D=     \exp \Bigl \{\int\limits _{w _P} ^{
w} A(z ,\tau )d\tau +
\int\limits _{z _P} ^{z} B(t,S(z)  )dt\Bigr \} \,\,
 \mbox{on }\widetilde{l}_2=\{{S}(z)=w _P\},
\end{eqnarray}
where  ${\hat L}^*$  is the adjoint operator  to $\hat L ^{\mathbb C}$ and
${\mathfrak R}(z_P ,w _P , z, w )$ is the Riemann function of
$\hat L$,

\begin{eqnarray}\nonumber
 {\hat L}^* {\mathfrak R}\equiv
\frac{\partial ^2 }{\partial z\partial w}{\mathfrak R}
-\pd{}{z}(A{\mathfrak R})
-\pd{}{w}(B{\mathfrak R}) +C{\mathfrak R}=0,\\
 {\mathfrak R}_ {|_{z=z_P}}= \exp \Bigl \{\int\limits _{w _P} ^{
w} A(z _P,\tau )d\tau \Bigr \}, \\
 {\mathfrak R}_ {|_{w =w _P}}= \exp \Bigl \{\int\limits _{z _P} ^{
z} B(t,w _P )dt \Bigr \}.
\end{eqnarray}
Note that for the Helmholtz equation 
$$
\Delta  _{x,y}u+\lambda ^2u=0,
$$
where $c=\lambda^2$ is a positive constant, $c_0(P,\Gamma )=1$, $a=b=0$,  and
${\mathfrak R}=J_0(\lambda\sqrt{(x-x_0)^2+(y-y_0)^2}\,)$. Thus, formula (\ref{RF}) could be slightly simplified, but it is still nonlocal \cite{nonloc}.

In the case of nonhomogeneous Dirichlet condition, when $u(x,y)=\phi (x,y)$ on $\Gamma$ with the same assumptions about $\phi$ as above, the reflection formula has the form \cite{ss}, p.471:
\begin{eqnarray}\label{RFHD}
\, & u(x_0,y_0) \,   =\,  -\,   u(R(x_0,y_0))+\,  \qquad\qquad\qquad\qquad\qquad\qquad\qquad\\
& +\frac{1}{2i}\int\limits   _{\Gamma }^{R(x_0,y_0)} \Bigl (
\bigl
\{
u\pd{V^D}{x}-V^D\pd{u}{x}\bigr \} dy
      - \bigl \{u\pd{V^D}{y}-V^D\pd{u}{y} \bigr \} dx
\Bigr )+\mathbb{F}^D[\phi (x_0,y_0)] .\nonumber
\end{eqnarray}
Here 
\begin{equation}
\mathbb{F}^D[\phi (x_0,y_0)]=\int\limits   _{\gamma _{\mathbb C}} \phi\Bigl (
\pd{}{x}(G- G^D)\, dy
      - \pd{}{y}(G- G^D)\, dx
\Bigr ),
\end{equation}
where $G=-\frac{1}{4\pi}(G_1+G_2)$ and 
$G^D=-\frac{1}{4\pi}( G_1^D+ G_2^D)$ are a fundamental solution and the reflected fundamental solution to the Helmholtz equation respectively, contour 
$\gamma _{\mathbb C}\subset \Gamma _{\mathbb C}$ that surrounds  both points
$(\widetilde{S}(w_0),w_0)$ and $(z_0, S(z_0))$.
Here
\begin{eqnarray}\label{fundH1}
G_1 & =\sum\limits _{k=0}^{\infty} \frac{[-\lambda
^2(z-z_0)(w -w _0)]^k}{4^k (k!)^2}\Bigl ( \ln (z-z_0 ) -C_k\Bigr ), \\
G_2 & =\sum\limits _{k=0}^{\infty} \frac{[-\lambda
^2(z-z_0)(w -w _0)]^k}{4^k (k!)^2}\Bigl ( \ln (w -w
_0 ) -C_k\Bigr ), \label{fundH2}
\end{eqnarray}
\begin{eqnarray}
 G_1^D & =
\sum\limits _{k=0}^{\infty} a_k^1 (z,w)\,
\frac{(\widetilde S(w)-z_0 )^k}{k!}\,
\Bigl ( \ln (\widetilde S(w)-z_0 ) -C_k\Bigr ), \\
 G_2^D & =\sum\limits _{k=0}^{\infty} 
a_k^2 (z,w)\,
\frac{(S(z)-w_0 )^k}{k!}\,
\Bigl ( \ln (S(z) -w
_0 ) -C_k\Bigr ), \label{fundHR}
\end{eqnarray}
where
$$
 C_0=0,\qquad C_k=\sum\limits _{l=1}^k \frac{1}{l},\qquad
 a^1_0=a^2_0=1,
$$
$$
\pd{a_{k+1}^1}{z}{\widetilde S^{\prime}(w)}=-\frac{\partial ^2 a_k^1}{\partial z\partial w}-
\frac{\lambda ^2}{4} a_k^1 , \quad
\pd{a_{k+1}^2}{w}S^{\prime}(z)=-\frac{\partial ^2 a_k^2}{\partial z\partial w}-
\frac{\lambda ^2}{4} a_k^2 ,
$$
$$
{a_k^1}_{|_{\Gamma _{\mathbb C}}}=\frac{[-\lambda
^2(w -w _0)]^k}{4^k k!}, \quad
{a_k^2}_{|_ {\Gamma _{\mathbb C}}}=\frac{[-\lambda
^2(z-z_0)]^k}{4^k k!}, \quad k=1,2,...\, .
$$
Introducing a notation $\omega (\cdot )=\pd{}{y}dx-\pd{}{x}dy$, formula (\ref{RFHD}) can be rewritten as
\begin{equation}
u(x_0,y_0)  =  - u(R(x_0,y_0))  +
\frac{1}{2i}\int\limits   _{\Gamma }^{R(x_0,y_0)} 
V^D\omega (u)-u\omega (V^D)+\int\limits _{\gamma _{\mathbb C}}\phi\omega(\widetilde G^D-G).
\end{equation}

\subsection{Reflection for nonhomogeneous Neumann data}

Homogeneous reflection formula for Helmholtz subject to Neumann condition 
can be written as a  following integro-differential operator \cite{S99},
\begin{equation}\label{E:1.2}
u(x_0 ,y_0 )=u(R(x_0 ,y_0 ))
\end{equation}
$$+\frac{1}{2i}\int\limits_{\Gamma }^{R(x_0 ,y_0 )}V^N(x,y,x_{0},y_{0})\omega(u(x,y))-u(x,y)
\omega(V^N(x,y,x_{0},y_{0})),$$ for any point $(x_0, y_0)$
sufficiently close to $\Gamma$. 
The integral is independent on the path
joining an arbitrary point on $\Gamma$ with the point $R(x_0 ,y_0)$ as above,
$$
V^N=V^N(x_P ,y_P ,x,y)=V_1^N (x_P ,y_P ,x,y)-V_2^N (x_P ,y_P ,x,y).
$$
Here functions $V_j^N$ are defined as the solutions to the following problems:
\begin{eqnarray}
  \frac{\partial ^2 V^N_j}{\partial z\partial w}+\frac{\lambda ^2}{4} V_j^N=0, \qquad j=1,2\, ,
\label{444}\\
 \omega ^*({V_j^N})_ {|_{\Gamma _{\mathbb{C}}}}=\omega ^*({\mathfrak R})_{|_{\Gamma _{\mathbb{C}}}},  \qquad  j=1,2\, ,\label{445}\\
 V_1^N= -1 \quad
 \mbox{on } \quad \widetilde{l}_1=\{\widetilde{S}(w )=z_0\},\\
 V_2^N=  -1   \quad
 \mbox{on } \quad  \widetilde{l}_2=\{{S}(z)=w _0\},\label{447}
\end{eqnarray}
where $\omega ^*(\cdot )=i(\pd{}{z}dz-\pd{}{w}dw)$.

Analogously to the Dirichlet case, a nonhomogeneous reflection formula has an extra term \cite{S99},
\begin{eqnarray}\label{E:1.2}
 u(x_0 ,y_0 )  =u(R(x_0 ,y_0 ))\qquad\qquad\qquad\qquad\qquad\qquad\qquad\\
+  \frac{1}{2i}\int\limits_{\Gamma }^{R(x_0 ,y_0 )}V^N\omega(u(x,y))-u(x,y)
\omega(V^N)+\mathbb{F}^N[\psi (x_0,y_0)],\nonumber
\end{eqnarray}
where
\begin{equation}
\mathbb{F}^N[\psi (x_0,y_0)]=\int\limits   _{\gamma _{\mathbb C} } \psi (z,S(z))( G^N-G)
\sqrt{S^{\prime}(z)}\, dz.
\end{equation}
Here 
$ G^N=-\frac{1}{4\pi}( G_1^N+ G_2^N)$, where
\begin{eqnarray}\label{fundHRN}
 G_1^N & =
\sum\limits _{k=0}^{\infty} b_k^1 (z,w)\,
\frac{(\widetilde S(w)-z_0 )^k}{k!}\,
\Bigl ( \ln (\widetilde S(w)-z_0 ) -C_k\Bigr ), \\
 G_2^N & =\sum\limits _{k=0}^{\infty} 
b_k^2 (z,w)\,
\frac{(S(z)-w_0 )^k}{k!}\,
\Bigl ( \ln (S(z) -w
_0 ) -C_k\Bigr ), 
\end{eqnarray}
$ b^1_0=b^2_0=-1$, for $k=1,2,...\,$,
$$
\pd{b_{k+1}^1}{z}{\widetilde S}^{\prime}(w)=-\frac{\partial ^2 b_k^1}{\partial z\partial w}-
\frac{\lambda ^2}{4} b_k^1 , \quad
\pd{b_{k+1}^2}{w}S^{\prime}(z)=-\frac{\partial ^2 b_k^2}{\partial z\partial w}-
\frac{\lambda ^2}{4} b_k^2 ,
$$
$$
{b_{k+1}^1}{\widetilde S}^{\prime}(w)_{|_ {\Gamma _{\mathbb C}}}={\widetilde S}^{\prime}(w) \pd{b_k^1}{z}-\pd{b_k^1}{w} _{|_ {\Gamma _{\mathbb C}}}
+\frac{(-\lambda ^2)^k(w-w_0)^{k-1}}{4^{k+1} (k+1)!}
(4k(k+1)+{\widetilde S}^{\prime}(w)\lambda ^2(w-w_0)^2)\, , 
$$
$$
{b_{k+1}^2} S^{\prime}(z)_{|_ {\Gamma _{\mathbb C}}}=S^{\prime}(z) \pd{b_k^2}{w}-\pd{b_k^2}{z} _{|_ {\Gamma _{\mathbb C}}}
+\frac{(-\lambda ^2)^k(z-z_0)^{k-1}}{4^{k+1} (k+1)!}
(4k(k+1)+S^{\prime}(z)\lambda ^2(z-z_0)^2) \, .
$$

\subsection{Representation of solutions to Cauchy problem for Helmholtz equation}

In this subsection we consider a Cauchy problem for Helmholtz equation,
\begin{eqnarray}\label{5_1ch}
\frac{\partial ^2 u}{\partial z\partial w}+\frac{\lambda ^2}{4} u=0\quad \mbox{near} \quad \Gamma_{\mathbb{C}},\\
u(z,w)=\phi (z,S(z))\quad \mbox{on} \quad \Gamma_{\mathbb{C}} ,\label{5_2ch}\\
\pd{u}{z}-\pd{u}{w}S^{\prime}=i \psi (z, S(z)){\sqrt{ S^{\prime} }}\quad \mbox{on} \quad \Gamma_{\mathbb{C}},\label{3ch}
\end{eqnarray}
 and discuss relationship between a representation of its solution with nonhomogeneous reflections
 mentioned above. 
Proceeding the same way as for the Laplace's equation, that is, eliminating $u(R(x_0 ,y_0 ))$
by adding formulae (\ref{RFHD}) and (\ref{E:1.2}), we obtain the following formula
\begin{eqnarray}\label{E:1.2C}\nonumber
u(x_0 ,y_0 )=
\frac{1}{4i}\int\limits_{\Gamma }^{R(x_0 ,y_0 )}(V^D+V^N)\omega(u(x,y))-u(x,y)
\omega(V^D+V^N)\\
+\frac{1}{2}\int\limits   _{\gamma _{\mathbb C}} \psi (z,S(z))( G^N-G)
\sqrt{ S^{\prime}(z)}\, dz +\phi\omega( G^D-G).\label{CHelp}
\end{eqnarray}
Unlike (\ref{CR}) this formula involves values of function $u$ on both sides of the curve $\Gamma$, so the natural question is whether or not the first integral vanishes.

Note that this formula reduces to (\ref{CR}) if $\lambda =0$. Indeed, in this case $V^D=V^N=0$, making the first integral vanish. Thus, formula (\ref{CHelp}) can be rewritten as
\begin{eqnarray}\nonumber
u(z_0 ,w_0 )=\frac{1}{2}\int\limits   _{\gamma _{\mathbb C}}\phi\omega^*( G^D_L-G_L)
\qquad\qquad\qquad\\
+\frac{1}{2}\int\limits   _{\gamma _{\mathbb C}} \psi (z,S(z))( G^N_L-G_L)
\sqrt{S^{\prime}(z)}\, dz ,\label{CHelpp}
\end{eqnarray}
where
$G_L=-\frac{1}{4\pi}\ln [(z-z_0)(w-w_0)]$,
$G_L^D=-\frac{1}{4\pi}\ln [(\widetilde S(w)-z_0)(S(z)-w_0)]$, and $G_L^N=-G_L^D$. Computation of  residues in the first integral and application of  properties of logarithms results in (\ref{CR}).

Another easy case to check is when $\Gamma$ is a line, 
 $\Gamma := \{ \alpha x+\beta y + \delta=0 \}$. Then the Schwarz function is $S(z)=mz+q$, 
where 
$$
m=\frac{\beta ^2 -\alpha ^2+i\, 2\alpha\beta}{\alpha ^2 +\beta ^2}, \qquad
q=\frac{ -2\alpha \delta +i\, 2\beta\delta}{\alpha ^2 +\beta ^2}.
$$
Functions $V_1^D$ and $V_2^D$ are equal,
\begin{equation}\nonumber
V_1^D=V_2^D =\sum\limits _{k=0}^{\infty}
 \frac{\bigl ( (mz+q-w _0)({\bar m}w +{\bar q} -z_0) (-\lambda ^2)\bigr )^k}{4^k(k!)^2}\,,
\end{equation}
and
therefore, $V^D=V_1^D-V_2^D=0$. Similarly $V_j^N=-V_j^D$, $j=1,2$ resulting in $V^N=0$. Thus, the first integral vanishes  as well. The second integral can be modified using the fact that $G^N=-G$ and $\omega ^* (G^D)=-\omega ^* (G)$ on $\Gamma _\mathbb{C}$,
which implies
\begin{equation}\label{prelim}\nonumber
u(x_0,y_0)=-\int\limits _{\gamma_{\mathbb C}}
G\psi\,\sqrt{S^{\prime}}\,dz+\phi\,\omega ^*(G) .
\end{equation}
It turns out that the latter formula holds for any algebraic curve $\Gamma$.
\begin{theorem}\label{TH}
Solution to the Cauchy problem (\ref{5_1ch})-(\ref{3ch})
is given by the formula
\begin{eqnarray}\label{CRH}
 u(z_0,w_0)=\frac{ \phi (\widetilde{S}(w_0),w_0)
+\phi (z_0, S(z_0))}{2}+\frac{i}{2} \int\limits _{ \widetilde{S}(w_0)}^{z_0}\phi(z, S(z))\omega ^*(J_0)\\
+ \frac{i}{2} \int\limits _{ \widetilde{S}(w_0)}^{z_0}J_0(\lambda \sqrt{(z-z_0)(S(z)-w_0)})
 \psi(z, S(z)) \sqrt{S^{\prime}}\, dz.
\nonumber
\end{eqnarray}
\end{theorem}
{\it Proof of Theorem \ref{TH}. }
We start with the Green's formula,
$$
u(x_0,y_0)=\int\limits _{\gamma}\Bigl ( 
u\pd{G}{n}-G\pd{u}{n} \Bigr )d\tau ,
$$ 
where $\gamma \subset\mathbb{R}^2$ is a small contour surrounding the point $(x_0,y_0)\in \mathbb{R}^2$ located close to $\Gamma$ and $G$ is a fundamental solution to the Helmholtz equation.
After complexification and deformation of the contour $\gamma$ to $\gamma _{\mathbb C}\subset\Gamma_{\mathbb C}$, we have
\begin{equation}\label{deformed}
u(x_0,y_0)=\int\limits _{\gamma_{\mathbb C}}
G\,\omega ^* (u)-u\,\omega ^*(G) .
\end{equation}
Taking into account conditions (\ref{5_2ch}) and (\ref{3ch}), the latter formula becomes
\begin{equation}\label{prelim}
u(x_0,y_0)=-\int\limits _{\gamma_{\mathbb C}}
G\psi\,\sqrt{S^{\prime}}\,dz+\phi\,\omega ^*(G) .
\end{equation}
Note that if we plug for $G$ a fundamental solution to the Laplace's equation,
$$
G_L=-\frac{1}{4\pi}\ln [(x-x_0)^2+(y-y_0)^2],
$$
and simplify the expression computing residues and using the properties of logarithms,
formula (\ref{prelim}) reduces to formula (\ref{CR}).
$$
u(z_0,w_0)=\frac{ \phi (\widetilde{S}(w_0),w_0)+\phi (z_0, S(z_0))}{2}+
\frac{i}{2} \int\limits _{ \widetilde{S}(w_0)}^{z_0}\psi(z, S(z)) \sqrt{S^{\prime}}\, dz.
$$
If we plug for $G$ the expression in terms of (\ref{fundH1}) and (\ref{fundH2}), 
formula (\ref{prelim}) can be rewritten as (\ref{CRH}).


\section{Elliptic growth}

\subsection{Solution to elliptic growth: Helmholtz equation}

Formula (\ref{CRH}) may be applied for solving elliptic growth problems, for instance, for the homogeneous screening mentioned in \cite{KhMinPut}.
Consider an elliptic growth problem involving the Helmholtz equation,
\begin{eqnarray}\label{2_1E}
\Delta p(x,y,t)+\lambda^2p(x,y,t) =\mu ,\\
p(x,y,t)=0 \quad \mbox{on} \quad \Gamma (t), \label{2_2E}\\
-k\pd{p}{n} =v_n \quad \mbox{on} \quad \Gamma(t),\label{3E}
\end{eqnarray}
where the support of the distribution $\mu$ is not near the free boundary. Setting $\phi =0$
and $\psi =i\dot S/(2k\sqrt{S^{\prime}})$, we obtain a representation of solution to this elliptic growth problem
as
\begin{eqnarray}\nonumber
 p(z_0,w_0,t)=
 -\frac{1}{4k} \int\limits _{ \widetilde{S}(w_0)}^{z_0}\dot S\, J_0(\lambda \sqrt{(z-z_0)(S(z)-w_0)})
 \, dz.
\end{eqnarray}

\subsection{Solution to elliptic growth: general case}

The above consideration implies that for a general elliptic equation of the second order 
(\ref{geneq}),  formula (\ref{deformed}) is replaced with
\begin{equation}
u(z_0,w_0)=\int\limits _{\gamma_{\mathbb C}}
G\,\omega ^* (u)-u\,\omega ^*(G) +2iGu(B\,dz-A\,dw),
\end{equation}
where $G$ is a fundamental solution of (\ref{geneq}).

Thus, the solution to the corresponding  Cauchy problem
\begin{eqnarray}\nonumber
\frac{\partial ^2 u}{\partial z\partial w}+A\pd{u}{z}+B\pd{u}{w} +Cu=0\quad \mbox{near} \quad \Gamma_{\mathbb{C}},\\
u(z,w)=\phi (z,S(z))\quad \mbox{on} \quad \Gamma_{\mathbb{C}} ,\nonumber\\
\pd{u}{z}-\pd{u}{w}S^{\prime}=i \psi (z, S(z)){\sqrt{ S^{\prime} }}\quad \mbox{on} \quad \Gamma_{\mathbb{C}},\nonumber
\end{eqnarray}
has the following representation
\begin{equation}
u(z_0,w_0)=-\int\limits _{\gamma_{\mathbb C}}
G\psi \sqrt{S^{\prime}}\, dz+\phi\,\omega ^*(G) +2i\phi G(B\,dz-A\,dw).
\end{equation}
Therefore, a solution to the elliptic growth problem
\begin{eqnarray}\label{2_1E}
\Delta p+a(x,y)\pd{p}{x}+b(x,y)\pd{p}{y}+c(x,y)p =\mu ,\\
p=0 \quad \mbox{on} \quad \Gamma (t), \label{2_2E}\\
-k\pd{p}{n} =v_n \quad \mbox{on} \quad \Gamma(t),\label{3E}
\end{eqnarray}
can be written as
$$
p(z_0,w_0,t)=
 -\frac{1}{4k} \int\limits _{ \widetilde{S}(w_0)}^{z_0}\dot S\, {\mathfrak R}(z_0 ,w _0 , z, S(z) )
 \, dz,
$$
where ${\mathfrak R}(z_0,w_0,z,w)$ is the Riemann function for (\ref{geneq}).

\section{Conclusions}

In this paper we discussed  the nonhomogeneous reflection formulae subject to the Dirichlet and Neumann conditions. We also showed their connections with continuation to the Cauchy's problems. The latter was applied to the Laplacian and the elliptic growth.

 

\providecommand{\bysame}{\leavevmode\hbox to3em{\hrulefill}\thinspace}

\end{document}